\documentclass[10pt]{article}
\usepackage{amsmath,amsfonts,amscd,amssymb}

\newtheorem{theorem}{Theorem}[section]
\newtheorem{lemma}[theorem]{Lemma}

\newtheorem{example}[theorem]{Example}
\newtheorem{remark}[theorem]{Remark}
\newtheorem{corollary}[theorem]{Corollary}

\newtheorem{assumption}[theorem]{Assumption}
{\catcode`@=11\global\let\c@equation=\c@theorem}

\newcommand{\zz}{{\mathbb Z}}

\newcommand{\cc}{{\mathbb C}}

\newcommand{\qq}{{\mathbb Q}}
\newcommand{\pp}{{\mathbb P}}
\newcommand{\rr}{{\mathbb R}}

\newcommand{\calb}{{\mathcal B}}
\newcommand{\call}{{\mathcal L}}
\newcommand{\caln}{{\mathcal N}}

\newcommand{\ch}{\operatorname{ch}}
\newcommand{\cone}{\operatorname{cone}}
\newcommand{\id}{\operatorname{id}}
\newcommand{\induc}{\operatorname{induc.}}
\newcommand{\im}{\operatorname{im}}
\newcommand{\pr}{\operatorname{pr}}
\newcommand{\sign}{\operatorname{sign}}

\newcommand{\squarematrix}[4]{\left( \begin{array}{cc} #1 & #2 \\ #3 &
#4
\end{array} \right)}

\providecommand{\proof}{\noindent{\bf \underline{\underline{Proof}} : }}
\providecommand{\qed}{\qquad \rule{2mm}{2mm} \bigskip}

\newcommand{\comment}[1]                      %comment of the author
{
{{\bf Comment: } {\ttfamily #1}}
}

\newcommand{\tit}[1]{\begin{bf} \begin{center} \begin{Large}
\section{#1}
\label{#1}
\end{Large}\end{center}\end{bf}
\nopagebreak}

\begin{document}

%\date{Last edited: March 24, 2000 or later\\last compiled: \today}

\title{On the cut and paste property of higher signatures of a closed
oriented manifold\\
by \\Eric Leichtnam and Wolfgang L\"uck\\
with an appendix\\
Mapping tori of special diffeomorphisms\\
by\\
Matthias Kreck}
\author{
Eric Leichtnam\thanks{\noindent e-mail:
leicht@math.jussieu.fr}\\ Institut de Chevaleret \\Plateau E, 7$^{eme}$
\'etage\\ (Alg\`ebres d' op\'erateurs)\\ 175 Rue de Chevaleret
\\75013  Paris\\France
\and Wolfgang L\"uck\thanks{\noindent email:
lueck@math.uni-muenster.de\protect\\
www: ~http://www.math.uni-muenster.de/u/lueck/org/staff/lueck/\protect}\\
Mathematisches Institut\\ Universit\"at M\"unster\\
Einsteinstr.~62\\ 48149 M\"unster\\Germany
\and Matthias Kreck
\thanks{\noindent email: kreck@mathi.uni-heidelberg.de \protect}\\
Mathematisches Institut \\ Universit\"at Heidelberg\\
Im  Neuenheimer Feld 288\\ 69120 Heidelberg\\Germany}

\maketitle
%%%%%%%%%%%%%%%%%%%%%%%%%%%% Abstract  %%%%%%%%%%%%%%%%%%%%%%%%%%%%%%%%%%%%%%%
\typeout{-----------------------  Abstract  ------------------------}
\begin{abstract}
We extend the notion of the symmetric
signature $\sigma(\overline{M},r) \in L^n(R)$
for a compact $n$-dimensional manifold $M$
without boundary, a reference map  $r:M \to BG$
and a homomorphism of rings with involutions
$\beta: \zz G \to R$ to the case with boundary
$\partial M$, where $(\overline{M},\overline{\partial M}) \to (M,\partial M)$ 
is the $G$-covering associated to $r$.
We need the assumption that
$C_*(\overline{\partial M}) \otimes_{\zz G} R$ is
$R$-chain homotopy equivalent to a $R$-chain
complex $D_*$ with trivial $m$-th differential for
$n = 2m$ resp. $n = 2m+1$. We prove a glueing formula, homotopy invariance
and additivity for this new notion. 

Let $Z$  be  a closed oriented
manifold with reference map $Z \to BG$. 
Let $F \subset Z$ be a cutting codimension one submanifold
$F \subset Z$ and let $\overline{F} \to F$ be the associated $G$-covering.
Denote by $\alpha_m(\overline{F})$
the $m$-th Novikov-Shubin invariant and by $b_m^{(2)}(\overline{F})$
the $m$-th $L^2$-Betti number. If for  the discrete group $G$ 
the Baum-Connes assembly map is rationally injective, then   
we use  $\sigma(\overline{M},r)$
to prove the additivity (or cut and paste property)
of the higher signatures of $Z$,
if we have  $\alpha_m(\overline{F}) = \infty^+$
in the case $n = 2m$ and, in the case $n = 2m+1$, 
if we have $\alpha_m(\overline{F}) = \infty^+$ and
$b_m^{(2)}(\overline{F}) = 0$. This additivity result 
had been proved (by a different method) in 
\cite[Corollary 0.4]{Leichtnam-Lott-Piazza(1999)}
when $G$ is  Gromov hyperbolic or virtually nilpotent.
We give new examples, where these conditions 
are not satisfied and additivity fails.

We explain at the end of the introduction why our paper is greatly motivated by
and partially extends some of the work of Leichtnam-Lott-Piazza 
\cite{Leichtnam-Lott-Piazza(1999)},
Lott \cite{Lott (1999)} and Weinberger \cite{Weinberger(1997)}.

\smallskip

\end{abstract}

%%%%%%%%%%%%%%%%%%%%%%%%%%%%%Introduction %%%%%%%%%%%%%%%%%%%%%%%%%%%%%%%%%%%%%%%
\typeout{-----------------------  Introduction ------------------------}

\setcounter{section}{-1}
\tit{Introduction}

Let $M$ be an oriented compact $n$-dimensional manifold possibly with boundary.
Let $G$ be a (discrete) group and $r: M \to BG$ be a (continuous)  reference map
to its classifying space.
Fix an (associative) ring $R$ (with unit and) with involution and a homomorphism
$\beta: \zz G \to R$ of rings with involution.  Let
$\overline{\partial M} \to \partial M$ and
 $\overline{M} \to M$ be the $G$-coverings
associated to the maps $r|_{\partial M}: \partial M \to BG$ and
$r: M \to BG$. 
Following
\cite[Section 4.7]{Lott (1992)} and
\cite[Assumption 1 and Lemma 2.3]{Leichtnam-Lott-Piazza(1999)}, 
we make an assumption about
$(\partial M,r|_{\partial M})$. \\

\begin{assumption} \label{assumption}
Let $m$ be the integer for which either $n = 2m$ or $n = 2m+1$.
Let $C_*(\overline{\partial M})$ be the cellular $\zz G$-chain complex.
Then we assume that the $R$-chain complex
$C_*(\overline{\partial M}) \otimes_{\zz G} R$ is
$R$-chain homotopy equivalent to a $R$-chain complex $D_*$ whose
$m$-th differential $d_m: D_m \to D_{m-1}$ vanishes.
\end{assumption}

We will discuss this assumption later
(see Lemma \ref{lem: vanishing middle differential and NS-invariants}).

We first consider the easier and more satisfactory case $n = 2m$.
Under Assumption \ref{assumption} we will assign (in Section 
\ref{Computations in symmetric $L$-groups}) for $n = 2m$ to $(M,r)$ the element
\begin{eqnarray}
& \sigma(\overline{M},r)  \in L^{2m}(R),&  \label{sigma(M,r)}
\end{eqnarray}
which we will call the symmetric signature,
in the symmetric $L$-group $L^{2m}(R)$
(Here and in the sequel we are considering the projective version
and omit in the standard notation $L^{2m}_p(R)$ the index $p$)

This element  $\sigma(\overline{M},r) $ agrees with the 
symmetric signature in the sense of
\cite[Proposition 2.1]{Ranicki(1980b)},
\cite[page 26]{Ranicki(1981)}, provided that $\partial M$ is empty.
If $\partial M$ is nonempty and 
$D_m = 0$ then $\sigma(\overline{M},r) $ was previously 
considered in \cite{Weinberger(1997)} and 
\cite[Appendix A]{Lott (1999)}.

The main properties of this invariant will be that it occurs in
a glueing formula, is a homotopy invariant
and is related to higher signatures as explained in Theorem
\ref{the: main properties of sigma(M,r)}, Theorem
\ref{the: main properties of sigma(M,r,L)} and
and Corollary
\ref{cor: Additivity of higher signatures}.

\begin{theorem} \label{the: main properties of sigma(M,r)}
\begin{enumerate}

\item Glueing formula
\label{the: main properties of sigma(M,r): glueing formula}
\\[1mm]
Let $M$ and $N$ be two oriented compact $2m$-dimensional
manifolds with boundary and let
$\phi: \partial M \to \partial N$ be an orientation
preserving diffeomorphism.
Let $r: M\cup_{\phi} N^- \to BG$ be a reference map.
Suppose that $(\partial M,r|_{\partial M})$ satisfies
Assumption \ref{assumption}. Then
$$\sigma(\overline{M\cup_{\phi}N^-},r) ~ = ~
\sigma(\overline{M},r|_M) - \sigma(\overline{N},r|_N);$$

\item Additivity \label{the: main properties of sigma(M,r): additivity}
\\[1mm]
Let $M$ and $N$ be two oriented compact $2m$-dimensional manifolds with boundary
and let
$\phi,\psi: \partial M \to \partial N$
be orientation preserving diffeomorphisms.
Let $r: M\cup_{\phi} N^- \to BG$ and $s: M\cup_{\psi} N^- \to BG$ be reference
maps such that $r|_M \simeq s|_M$ and $r|_N \simeq s|_N$ holds, where $\simeq$
means homotopic. Suppose
that $(\partial M,r|_{\partial M})$ satisfies Assumption \ref{assumption}.
Then
$\sigma(\overline{M\cup_{\phi}N^-},r) ~ = ~ \sigma(\overline{M\cup_{\psi}N^-},s);$

\item Homotopy invariance
\label{the: main properties of sigma(M,r): homotopy invariance}
\\[1mm]
Let $M_0$ and $M_1$ be two oriented compact $2m$-dimensional manifolds
possibly with boundaries together with
reference maps $r_i: M_i \to BG$ for $i = 0,1$.
Let $(f,\partial f): (M_0,\partial M_0) \to (M_1,\partial M_1)$
be an orientation preserving homotopy equivalence
of pairs with $r_1 \circ f \simeq r_0$. Suppose that
$(\partial M_0,r_0|_{\partial M_0})$ satisfies
Assumption \ref{assumption}. Then
$$\sigma(\overline{M_0},r_0) = \sigma(\overline{M_1},r_1).$$

\end{enumerate}
\end{theorem}

Next we consider the case $n= 2m+1$. Then we need besides
Assumption \ref{assumption} the following additional input.
Assumption \ref{assumption} implies that
$H_m(C_*(\overline{\partial M}) \otimes_{\zz G} R)$ is a finitely generated
projective $R$-module and that we get from Poincar\'e duality the structure
of a (non-degenerate) $(-1)^m$-symmetric form $\mu$ on it.
Following
\cite[Section 3]{Leichtnam-Lott-Piazza(1999)}, 
we will assume that we have specified a stable Lagrangian
$L \subset H_m(C_*(\overline{\partial M}) \otimes_{\zz G} R)$.
The existence of a stable Lagrangian follows automatically
if $2$ is a unit in $R$ (see Lemma
\ref{lem: existence of stable Lagrangian}).
Under Assumption \ref{assumption}
and after the  choice of a stable Lagrangian $L$
we can assign for $n = 2m+1 $ to
$(M,r,L)$ an element, which we will call the symmetric signature,
 in the symmetric $L$-group $L^{2m+1}(R)$ (see Section 
\ref{Computations in symmetric $L$-groups})
\begin{eqnarray}
& \sigma(\overline{M},r,L)  \in L^{2m+1}(R).&  \label{sigma(M,r,L)}
\end{eqnarray}
It agrees with the symmetric signature in the sense of
\cite[Proposition 2.1]{Ranicki(1980b)},
\cite[page 26]{Ranicki(1981)}, provided that $\partial M$ is empty.

\begin{theorem} \label{the: main properties of sigma(M,r,L)}

\begin{enumerate}

\item Glueing formula
\label{the: main properties of sigma(M,r,L): glueing formula}
\\[1mm]
Let $M$ and $N$ be oriented compact $(2m+1)$-dimensional manifolds with boundary and let
$\phi: \partial M \to \partial N$ be an orientation
preserving diffeomorphism.
Let $r: M\cup_{\phi} N^- \to BG$ be a reference map.
Suppose that $(\partial M,r|_{\partial M})$ satisfies
Assumption \ref{assumption}. Suppose that we have fixed two stable Lagrangians
$K \subset H_m(C_*(\overline{\partial M}) \otimes_{\zz G} R)$
and $L \subset H_m(C_*(\overline{\partial N}) \otimes_{\zz G} R)$
such that the isomorphism $H_m(C_*(\overline{\partial M}) \otimes_{\zz G} R)
\xrightarrow{\cong} H_m(C_*(\overline{\partial N}) \otimes_{\zz G} R)$
of $(-1)^m$-symmetric forms induced by $\phi$ sends  $K$ to $L$ stably.
Then
$$\sigma(\overline{M\cup_{\phi}N^-},r) ~ = ~
\sigma(\overline{M},r|_M,K) - \sigma(\overline{N},r|_N,L);$$

\item Additivity
\label{the: main properties of sigma(M,r,L): additivity}
\\[1mm]
Let $M$ and $N$ be oriented compact $(2m+1)$-dimensional manifolds with
boundary and let
$\phi,\psi: \partial M \to \partial N$ be two orientation preserving diffeomorphisms.
Let $r: M\cup_{\phi} N^- \to BG$ and $s: M\cup_{\psi} N^- \to BG$ be reference
maps together with homotopies $h_M: r|_M \simeq s|_M$ and 
$h_M: r|_N \simeq s|_N$. Suppose
that $(\partial M,r|_{\partial M})$ satisfies Assumption \ref{assumption}.
Fix a stable Lagrangian
$K \subset H_m(C_*(\overline{\partial M}) \otimes_{\zz G} R)$.
The restriction of the homotopies $h_M$ and $h_N$ to $\partial M$ and $\partial N$ induce
a homotopy $k: r|_{\partial M} \circ \psi^{-1} \circ \phi
\simeq r|_{\partial M}$. We get from $\psi^{-1} \circ \phi$ and $k$ an
automorphism of the $(-1)^m$-symmetric form
$(H_m(C_*(\overline{\partial M}) \otimes_{\zz G} R),\mu)$. 
Let $L \subset H_m(C_*(\overline{\partial M}) \otimes_{\zz G} R)$
be the stable Lagrangian which is the image of $K$ under this automorphism. 
Thus we get a formation
$(H_m(C_*(\overline{\partial M}) \otimes_{\zz G} R),\mu,K,L)$
which defines an element in
$[H_m(C_*(\overline{\partial M}) \otimes_{\zz G} R),\mu,K,L]
\in L^{2m+1}(R)$ by suspension.
Then
$$\sigma(\overline{M\cup_{\phi}N^-},r) -
 \sigma(\overline{M\cup_{\psi}N^-},s) \,=
 \, [H_m(C_*(\overline{\partial M}) \otimes_{\zz G} R),\mu,K,L];
 $$
\item Homotopy invariance
\label{the: main properties of sigma(M,r,L): homotopy invariance}
\\[1mm]
Let $M_0$ and $M_1$ be oriented compact $(2m+1)$-dimensional manifolds
possibly with boundaries together with
reference maps $r_i: M_i \to BG$ for $i=0,1$. Let
$(f,\partial f): (M_0,\partial M_0) \to (M_1,\partial M_1)$ be an
orientation preserving  homotopy equivalence
of pairs together with a homotopy $h: r_1 \circ f \simeq r_0$. Suppose that
$(\partial M_0,r_0|_{\partial M_0})$ satisfies Assumption \ref{assumption}.
Suppose that we have fixed stable Lagrangians
$L_0 \subset H_m(C_*(\overline{\partial M_0}) \otimes_{\zz G} R)$ and
$L_1 \subset H_m(C_*(\overline{\partial M_1}) \otimes_{\zz G} R)$.
Let $L_0'$ be the image of $L_0$ under the isomorphism of
$(-1)^m$-symmetric forms
$(H_m(C_*(\overline{\partial M_0}) \otimes_{\zz G} R),\mu_0)
\xrightarrow{\cong}
(H_m(C_*(\overline{\partial M_1}) \otimes_{\zz G} R),\mu_1)$ induced by
$\partial f$ and the restriction of the homotopy $h$ to $\partial M_0$. 
We get a stable formation
$(H_m(C_*(\overline{\partial M_1}) \otimes_{\zz G} R),\mu_1, L_0',L_1)$
and thus by suspension an element

\noindent
$[H_m(C_*(\overline{\partial M_1}) \otimes_{\zz G} R),\mu_1, L_0',L_1]
\in L^{2m+1}(R)$. Then
$$\sigma(\overline{M_0},r_0,L_0) - \sigma(\overline{M_1},r_1,L_1)=
[H_m(C_*(\overline{\partial M_1}) \otimes_{\zz G} R),\mu_1, L_0',L_1].$$

\end{enumerate}
\end{theorem}

Of particular interest is the case, where $R$ is the real reduced
group $C^*$-algebra $C^*_r(G,\rr)$ or the complex 
reduced group $C^*$-algebra $C^*_r(G)$
and $\beta$ is the canonical map. Then
Assumption \ref{assumption} is equivalent to the assertion that the
$m$-th Novikov-Shubin invariant of $\overline{\partial M}$ is $\infty^+$
in the sense of \cite[Definition 1.8, 2.1 and 3.1]{Lott-Lueck (1995)}
(see Lemma \ref{lem: vanishing middle differential and NS-invariants}), and
the symmetric $L$-groups are $2$-periodic.
Moreover, the invariant $\sigma(\overline{M},r)$ is linked to higher
signatures as follows, provided that $\partial M$ is empty.

Recall that the higher signature $\sign_u(M,r)$ of a closed oriented
manifold $M$ with a reference map $r: M \to BG$
for a given class $u \in H^k(BG;\qq)$ is the rational number
$\langle \call(M) \cup r^* u,[M]\rangle$, where
$\call(M)\in \oplus_{i \ge 0} H^{4i}(M;\qq)$  is the
$L$-class of $M$, $[M]\in H_{\dim(M)}(M;\qq)$
is the homological fundamental class of $M$ and
$\langle ~,~ \rangle$ is the Kronecker pairing.
We will consider the following
commutative square of $\zz/4$-graded rational vector spaces
$$
\begin{CD}
\left(\Omega_*(BG)\otimes_{\Omega_*(*)} \qq\right)_n
@> \overline{D} > \cong >KO_n(BG)\otimes_{\zz} \qq @>\induc >>
K_n(BG) \otimes_{\zz} \qq
\\
@V \sigma VV @VA_{\rr}  VV @V A VV
\\
L^n(C_r^*(G;\rr))\otimes_{\zz} \qq @> \sign > \cong >
KO_n(C_r^*(G;\rr))\otimes_{\zz} \qq @> \induc >>
K_n(C_r^*(G)) \otimes_{\zz} \qq
\end{CD}
$$
Some explanations are in order.
We denote by $\qq$ the $\zz$-graded vector space which is
$\qq$ in each dimension divisible by four and zero elsewhere.
It can be viewed as a graded module over the $\zz$-graded ring
$\Omega_*(*)$ by the signature. Then the $\zz$-graded $\qq$-vector space
$\Omega_*(BG)\otimes_{\Omega_*(*)} \qq$ is four-periodic (by crossing with
$[\cc\pp^2]$) and hence
can be viewed as a $\zz/4$-graded vector space. The map
$\overline{D}$ is induced by the $\zz$-graded homomorphism
$$D: \Omega_n(BG) \to KO_n(BG),$$
which sends $[r: M \to BG]$ to the $K$-homology
class of the signature operator of the covering $\overline{M} \to M$
associated to $r$. The homological
Chern character is an isomorphism of $\zz/4$-graded rational vector spaces
$$\ch: KO_n(BG) \otimes_{\zz} \qq \xrightarrow{\cong}
\oplus_{k \ge 0} H_{4k+n}(BG;\qq).$$
By the Atiyah-Hirzebruch index theorem
the image $\ch \circ D([M,\id: M \to M])$
of the K-homology class of the signature operator of $M$ in $K_{\dim(M)}(M)$
under the homological Chern character $\ch$ is $\call(M) \cap [M]$.
This implies for any class $u \in H^k(BG;\qq)$
\begin{eqnarray}
\sign_u(M,r) & := & \langle \call(M) \cup r^* u,[M]\rangle \nonumber
\\ &  = & \langle r^* u,\call(M) \cap [M]\rangle \nonumber
\\ &  = & \langle u,r_*(\call(M) \cap [M]) \rangle \nonumber
\\ & = & \langle u,\ch \circ D([M,r])\rangle.
\label{index and higher signatures}
\end{eqnarray}
Hence the composition
$\ch \circ \overline{D}: \left(\Omega_*(BG)\otimes_{\Omega_*(*)} \qq\right)_n
\to \oplus_{k \ge 0} H_{4k+n}(BG;\qq)$
sends $[r: M \to BG] \otimes 1$ to the image
under $H_*(r): H_*(M;\qq) \to H_*(BG;\qq)$
of the Poincar\'e dual
$\call(M) \cap [M] \in \oplus_{i \ge 0}
H_{4i-\dim(M)}(M;\qq)$ of the $L$-class $\call(M)$. 

The map $\overline{D}$ is an isomorphism 
since it is a transformation of homology theories
\cite[Example 3.4]{Landweber(1976)} and induces an isomorphism for 
the space consisting of one point.
The map  $\sigma$ assigns to 
$[M,r]$ the associated symmetric Poincar\'e
$C_r^*(G;\rr)$-chain complex
$C_*(\overline{M}) \otimes_{\zz G} C_r^*(G;\rr)$.
The map $A_{\rr}$ resp. $A$ are assembly maps given by taking the index with
coefficients in $C^*_r(G;\rr)$ resp. $C^*_r(G)$.
The map $\sign$ is in dimension $n = 0 \mod 4$
given by taking the signature of a
non-degenerate symmetric bilinear form.
Notice that the map $\sign$ is  bijective by results of Karoubi
(see \cite[Theorem 1.11]{Rosenberg(1995)}).
The maps $\induc$ are given by induction with the inclusion
$\rr \to \cc$ and are injective. Obviously the right square commutes.
In order to show that the diagram commutes it suffices to prove this
for the outer square. Here the claim follows from the commutative diagram in
\cite[page 81]{Kasparov(1993)}.

The Novikov Conjecture says
that $\sign_u(M,r)$ is a homotopy invariant, i.e.
if $r: M \to BG$ and $s: N \to BG$ are closed orientable
manifolds with reference maps to $BG$ and $f: M \to N$
is a homotopy equivalence with $s \circ f \simeq r$, then
$\sign_u(M,r) = \sign_u(N,s)$.
Since the homological Chern character
is rationally an isomorphism for $CW$-complexes, one can say by
\eqref{index and higher signatures} that $D(M,r)$ is rationally
the same as the collection of all higher 
signatures. Moreover, the Novikov Conjecture
is equivalent to the statement that two elements $[M,r]$ and $[N,s]$
in $\Omega_n(BG)$ represent the same element in 
$\left(\Omega_*(BG)\otimes_{\Omega_*(*)} \qq\right)_n$ resp.
$KO_n(BG)\otimes_{\zz} \qq$ resp. $K_n(BG) \otimes_{\zz} \qq$, if they are homotopy 
equivalent. 

Notice that the Baum-Connes Conjecture for  $C_r^*(G)$
implies that $A$ and hence $A_{\rr}$ are rationally injective 
by the following argument (see
\cite[section 7]{Baum-Connes-Higson} for details).
The map $A$ can be written as the composition of the
map $K_n(BG) = K_n^G(EG) \to K_n(\underline{E}G)^G$,
which is given by the canonical map from
$EG$ to the classifying space $\underline{E}G$ of
proper $G$-actions and is always rationally
injective, and the Baum-Connes index map
$K^G_*(\underline{E}G) \to K_*(C^*_r(G))$,
which is predicted to be bijective by the Baum-Connes Conjecture.
Notice that $\sigma$ is injective if and only if $A_{\rr}$ is injective and that 
the injectivity of $A$ implies the injectivity of $A_{\rr}$ and hence of $\sigma$.
Since for a closed oriented manifold $M$ with reference map $r: M \to BG$
the image of $[r:M \to BG]$ under $\sigma$ (and
$\sign \circ\, \sigma $) is a homotopy invariant of
$r:M \to BG$, the commutativity of
the diagram above and  the rational injectivity of $A_{\rr}$
implies the homotopy invariance of \eqref{index and higher signatures}
and thus the Novikov conjecture.
Moreover, if $A_{\rr}$ is rationally injective, 
$D([M,r])$ contains rationally the same
information as $\sigma([M,r])$.
We mention that the Baum-Connes Conjecture and thus the rational injectivity of $A$
is known for a large class of groups, namely for all a-T-menable groups
\cite{Higson-Kasparov(1997)}. The rational injectivity of
$A$ is also known for all Gromov-hyperbolic groups \cite{Yu (1998)}.

From Theorem \ref{the: main properties of sigma(M,r)}
\eqref{the: main properties of sigma(M,r): additivity} and Theorem
\ref{the: main properties of sigma(M,r,L)}
\eqref{the: main properties of sigma(M,r,L): additivity}, we obtain
the following corollary which extends 
\cite[Corollary 0.4]{Leichtnam-Lott-Piazza(1999)} to more general groups $G$. 

\begin{corollary} \label{cor: Additivity of higher signatures}
Let $M$ and $N$ be two oriented compact $n$-dimensional manifolds with boundary
 and let
$\phi,\psi: \partial M \to \partial N$
be orientation preserving diffeomorphisms.
Let $r: M\cup_{\phi} N^- \to BG$ and $s: M\cup_{\psi} N^- \to BG$ be reference
maps such that $r|_M \simeq s|_M$ and $r|_N \simeq s|_N$ holds.
Denote by $\overline{\partial M} \to \partial M$ the $G$-covering associated to
$r|_{\partial M}: \partial M \to BG$.
If $n = 2m$, we assume for the $m$-th Novikov Shubin invariant
$\alpha_m(\overline{\partial M}) = \infty^+$ in the sense of
\cite{Lott-Lueck (1995)}. If $n = 2m+1$, we assume
$\alpha_m(\overline{\partial M}) = \infty^+$
and for the $m$-th $L^2$-Betti number $b_m^{(2)}(\overline{\partial M}) = 0$.
(We could replace the condition $b_m^{(2)}(\overline{\partial M}) = 0$
by the weaker but harder to check assumption that the automorphism in Theorem
\ref{the: main properties of sigma(M,r,L)}
\eqref{the: main properties of sigma(M,r,L): additivity}
induced by
$\psi^{-1}\circ \phi$ and the homotopy $k$ preserve (stably) a Lagrangian of
$H_m(\,C_*(\overline{\partial M})\otimes_{\zz G} C^*_r(G)\,) $.
Suppose furthermore that the map $A_{\rr}: KO_n(BG) \to K_n(C^*_r(G;\rr)$ is injective.
Then all the higher signatures are additive in the sense that
we have for all $u \in H^*(BG,\qq)$
\begin{eqnarray} \label{add}
\sign_u(M\cup_{\phi}N^-,r) = \sign_u(M\cup_{\psi}N^-,s).
\end{eqnarray}
\end{corollary}

We will construct in Example \ref{exa: non-additivity} many new examples
(especially for odd dimensional manifolds) of pairs of cut and paste 
manifolds $[ M\cup_{\phi} N^- ,r]$
 and $[M\cup_{\psi}N^-, s ]$ (even with $\partial M$ connected)
  such that $r|_M \simeq s|_M$ and $r|_N \simeq s|_N$ holds and for which 
   there exist higher signatures which do not 
satisfy (\ref{add}), ie  are not additive.
 There, the assumptions of 
 Corollary \ref{cor: Additivity of higher signatures},
are not fully  satisfied. 
The fact that, in general, higher signatures of closed manifolds are
not cut and paste invariant over $BG$ in the sense of 
\cite{Karras-Kreck-Neumann-Ossa (1973)},
was known before (see, for instance, \cite[Section 4.1]{Lott (1999)}).
The relationship to symmetric signatures of manifolds-with-boundary, and
to the necessity of Assumption \ref{assumption},
was pointed out by Weinberger (see  \cite[Section 4.1]{Lott (1999)}).
The problem was raised in \cite[Section 4.1]{Lott (1999)} 
of determining which higher 
signatures of closed manifolds 
are cut and paste invariant;  we refer to \cite[Section 4.1]{Lott (1999)} for
further discussion. It is conceivable that our Lemma 
\ref{lem: properties of sigma of a pair: Additivity and mapping tori} 
might help to provide, in the future, an answer to this problem.

Finally we explain why our paper is greatly motivated by
and related to the work of Leichtnam-Lott-Piazza \cite{Leichtnam-Lott-Piazza(1999)},
Lott \cite{Lott (1999)} and Weinberger \cite{Weinberger(1997)}.

The relevance of a gap condition in the middle degree on the boundary, when
considering topological questions concerning manifolds with boundary,
comes from Section 4.7 of Lott's paper \cite{Lott (1992)}.
This leads to Assumption 1 in the paper of
Leichtnam-Lott-Piazza \cite{Leichtnam-Lott-Piazza(1999)}
which, from \cite[Lemma 2.3]{Leichtnam-Lott-Piazza(1999)}, 
is virtually identical to our
Assumption \ref{assumption} and was the motivation for our 
Assumption \ref{assumption}.
Our construction of the invariant
$\sigma(\overline{M},r)$ by glueing algebraic Poincar\'e bordisms
is motivated by and extends the one of
Weinberger \cite{Weinberger(1997)}
(see also \cite[Appendix A]{Lott (1999)})
who uses the more restrictive assumption
that $C_*(\overline{\partial M}) \otimes_{\zz G} R$ is $R$-chain
homotopy equivalent to a $R$-chain complex $D_*$ with $D_m = 0$.
Notice that Weinberger's assumption
implies both our Assumption \ref{assumption} and
$H_m(C_*(\overline{\partial M}) \otimes_{\zz G} R) = 0$
so that there is only one choice of Lagrangian, namely $L = 0$.
The idea of using a Lagrangian subspace, instead of assuming the vanishing
of the relevant middle (co-)homology group, is taken from
Section 3 of \cite{Leichtnam-Lott-Piazza(1999)}.

In the case when $R = C^*_r(G)$ and under Assumption \ref{assumption}, 
an analog of our symmetric signature
$\sigma(\overline{M},r)$ was previously constructed
in \cite{Leichtnam-Lott-Piazza(1999)} as a conic index class
$\sigma_{conic} \in K_0( C^*_r(G))$. The homotopy invariance of
$\sigma_{conic}$, i.e. the analog of our Theorem 
\ref{the: main properties of sigma(M,r)}(c), was demonstrated in
\cite[Theorem 6.1]{Leichtnam-Lott-Piazza(1999)}.
Furthermore,
the analog of the right-hand-side of the equation in our Theorem
\ref{the: main properties of sigma(M,r,L)}(c) previously appeared
in \cite[Proposition 3.7]{Leichtnam-Lott-Piazza(1999)}

If $\calb^{\infty}$ is a smooth subalgebra of
$C^*_r(G)$, the authors of \cite{Leichtnam-Lott-Piazza(1999)}
computed the Chern character $\ch(\sigma_{conic}) \in 
\overline{H}_*(\calb^{\infty})$ of their conic index
explicitly in terms of the $L$-form of $M$
and a higher eta-form of $\partial M$. They identified
$\ch(\sigma_{conic})$ with the $\overline{H}_*(\calb^{\infty})$-valued
higher signature of $M$ introduced
in \cite{Lott (1992)}. 
From this the authors of
\cite{Leichtnam-Lott-Piazza(1999)} deduced 
their main result \cite[Theorem 0.1]{Leichtnam-Lott-Piazza(1999)}, namely
the homotopy invariance of the $\overline{H}_*(\calb^{\infty})$-valued
higher signature of a manifold with
boundary, as opposed to just the homotopy invariance of the
``symmetric signature'' $\sigma_{conic}$. (In fact, 
this was the motivation for the
use of $\sigma_{conic}$ in \cite{Leichtnam-Lott-Piazza(1999)},
instead of $\sigma(\overline{M},r)$.)

As an immediate consequence of its main result,
the paper \cite{Leichtnam-Lott-Piazza(1999)} deduced
the additivity of ordinary higher signatures of closed manifolds 
under Assumption \ref{assumption}, i.e. our Corollary
\ref{cor: Additivity of higher signatures}, 
in the case when $G$ is Gromov-hyperbolic or virtually nilpotent,
or more generally when 
$C^*_r(G)$ admits a smooth subalgebra $\calb^{\infty}$ with the property
that all of the group cohomology of $G$ extends to cyclic cocycles
on $\calb^{\infty}$
\cite[Corollary 0.4]{Leichtnam-Lott-Piazza(1999)}. It is known that
the Baum-Connes assembly map is rationally injective for such groups $G$.

What is new in our paper is the direct and purely algebraic
approach to the cut and paste
problem of higher signatures of closed manifolds,
through the construction of symmetric signatures
for manifolds with boundary under Assumption \ref{assumption}.
This leads to the main goal of the present paper,
namely, to the proof of the additivity of higher signatures 
under Assumption \ref{assumption} in
Corollary \ref{cor: Additivity of higher signatures}, provided
that the Baum-Connes assembly map is rationally injective.
In this way, our main result,
Corollary \ref{cor: Additivity of higher signatures}, is an extension of
\cite[Corollary 0.4]{Leichtnam-Lott-Piazza(1999)}.

%%%%%%%%%%%%%%%%%%%%%%%%%%%% Section 1 %%%%%%%%%%%%%%%%%%%%%%%%%%%%%%%%%%%%%%%
\typeout{-----------------------  Section 1 ------------------------}

\tit{Additivity and mapping tori in the bordism group}

Throughout this section $X$ is some topological space.
Denote by
$\Omega_n(X)$ the bordism group of closed oriented smooth $n$-dimensional
manifolds $M$ together with
a reference map $r: M \to X$. Consider quadruples $(F,h,r,H)$
consisting of a closed oriented $(n-1)$-dimensional manifold $F$ together with
an orientation preserving self-diffeomorphism $h: F \to F$, a reference
map $r: F \to X$  and a homotopy 
$H: F\times [0,1] \rightarrow X$ such that $H(-,0) = r$ and
$H(-,1) = r\circ h$.
The mapping torus $T_h$ is obtained from the cylinder $F \times [0,1]$
by identifying the bottom and the top by $h$, i.e. $(h(x),0) \sim (x,1)$.
This is again a closed smooth manifold and inherits a preferred orientation.
The map $r$  and the homotopy
$H$ yield a reference map $r_{T_h}: T_h \to X$
in the obvious way. Hence we can associate to such a quadruple an element
\begin{eqnarray}
[F,h,r,H] & := & [T_h,r_{T_h}] \hspace{3mm} \in \Omega_n(X).
\label{bordism element given by mapping torus}
\end{eqnarray}
Given two quadruples $(F,h[0],r[0],H[0])$
and $(F,h[1],r[1],H[1])$
with the same underlying manifold,
a homotopy between them is given by a family of such quadruples
$(F,h[t],r[t],H[t])$ for $t \in [0,1]$ such that
the family $h[t]: F \to F$ is a diffeotopy.
 One easily checks that for two    homotopic quadruples (as above) we have
in $\Omega_n(X)$
\begin{eqnarray}
[F,h[0],r[0],H[0]] & = &
[F,h[1],r[1],H[1]].
\label{homotopic quadruples give same bordism elements}
\end{eqnarray}
The required cobordism has as underlying manifold
$F\times [0,1] \times [0,1]/\sim$, where $\sim$
is the equivalence relation generated by
$(x,0,t) \sim (h^{-1}[t](x),1,t)$.

Given two quadruples of the shape  $(F,h,r,H)$ and
$(F,g,r,G)$, we can compose them to a quadruple
$(F,g \circ h, r, H\ast G)$, where $H\ast G$ is the
obvious homotopy $r \simeq r \circ g \circ h$
obtained from stacking together
$H$ and $G \times (h \times \id_{[0,1]})$.
One easily checks that  in $\Omega_n(X)$
\begin{eqnarray}
[F,g \circ h, r, H\ast G] & = &
[F,h,r,H] + [F,g, r, G].
\label{addition formula for quadruples}
\end{eqnarray}
The desired cobordism has as underlying
manifold $F \times [0,1] \times[0,1]/\sim$,
where $\sim$ is generated by $(x,0,t) \sim (h^{-1}(x),1,t)$ for
$t \in[0,1/3]$ and
$(x,1,t) \sim (g^{-1}(x),0,t)$ for $t \in[1/3,1]$.
We recognize the mapping torus of
$g \circ h$ as the part of the boundary
which is the image under the canonical projection
of the union of $F \times \{0\} \times[1/3,2/3]$,
$F \times [0,1] \times \{1/3\}$,
$F \times \{1\} \times [1/3,2/3]$ and  $F \times [0,1] \times \{2/3\}$.

Notice that the class of a quadruple in $\Omega_n(X)$
does depend on the choice of the homotopy.
Namely, consider two quadruples
$(F,h,r,H)$ and $(F,h,r,G)$
which differ only in the choice of the homotopy.
Let $u: F\times S^1 \to X$ be the obvious
 map induced by $r$ and composition of homotopies
$H \ast G^-: r \simeq r$. Then we get
from \eqref{homotopic quadruples give same bordism elements}
 and \eqref{addition formula for quadruples} in
$\Omega_n(X)$
\begin{eqnarray}
[F,h,r,H] - [F,h,r,G] & = & [u: F \times S^1 \to X].
\label{dependence of the bordism class of a quadruple on the homotopy}
\end{eqnarray}
The right side of
\eqref{dependence of the bordism class of a quadruple on the homotopy}
is not zero in general. Take for
instance $F = \cc\pp^2$ and $X = S^1$ and let
$u_k:F \times S^1 \to S^1$ be the
composition of the projection $F \times S^1 \to S^1$
with a map $S^1 \to S^1$ of degree $k \in \zz$.
Then in this situation  the right-handside of
 \ref{dependence of the bordism class of a quadruple on the homotopy}
becomes
 $[u_k: \cc\pp^2 \times S^1 \to S^1]$ for $r: \cc\pp^2 \rightarrow S^1$
a constant map. This element $[u_k: \cc\pp^2 \times S^1 \to S^1]$
is mapped under the isomorphism
$\Omega_5(S^1) \cong \Omega_4(*) \oplus \Omega_5(*) = \zz \oplus \zz/2$
to $(k,0)$.

Let $M$ and $N$ be compact oriented
$n$-dimensional manifolds and let $\phi,\psi: \partial M \to
\partial N$ be two orientation preserving diffeomorphisms.
By glueing we obtain closed oriented
$n$-dimensional manifolds $M \cup_{\phi} N^-$ and
$M \cup_{\psi} N^-$. Let $r: M \cup_{\phi} N^- \to X$
and $s: M \cup_{\phi} N^- \to X$
be two reference maps such that there exists two homotopies
 $H: r|_M \simeq s|_M$ and
$G: r|_N \simeq s|_N$.
By restriction we obtain homotopies
$H|_{\partial M}: r|_{\partial M} \simeq s|_{\partial M}$ and
$G|_{\partial N}: r|_{\partial N} \simeq s|_{\partial N}$. Notice that
(by construction)
$r|_{\partial N} \circ \phi = r|_{\partial M}$ and
$s|_{\partial N} \circ \psi = s|_{\partial M}$. Thus $H|_{\partial M}$ and
$G^-|_{\partial N} \circ (\psi \times \id)$
can be composed to a homotopy
$K: r|_{\partial M} \simeq r|_{\partial M} \circ \phi^{-1} \circ \psi$.
 Thus we obtain a quadruple
$(\partial M,\phi^{-1} \circ \psi, r|_{\partial M},K)$
in the sense of \eqref{bordism element given by mapping torus}.

\begin{lemma}
\label{lem: cutting and mapping torus formula in bordims group}
We get in $\Omega_n(X)$
\begin{eqnarray*}
[r: M \cup_{\phi} N^- \to X] - [s: M \cup_{\psi} N^- \to X] & = &
[\partial M,\phi^{-1} \circ \psi, r|_{\partial M},K].
\end{eqnarray*}
\end{lemma}
\proof
The underlying manifold  of the required bordism is
obtained by glueing parts of the boundary of
$M\times [0,1]$ and of  $N^-\times [0,1]$ together as described as follows.
Identify $(x,t) \in \partial M \times [0,1]$
with $(\phi(x),t)$ in $\partial N \times [0,1]$
if $0 \le t \le 1/3,$  and with $(\psi(x),t)$
in $\partial N$ if $2/3 \le t \le 1$. \qed

\begin{corollary} \label{cor: additivity in the bordism ring}
Suppose in the situation of Lemma 
\ref{lem: cutting and mapping torus formula in bordims group}
that $X = BG$ for a discrete group $G$  and that the image $H$ of the composition
$$\pi_1(\partial M) \to \pi_1(M) \xrightarrow{(r|_M)_*} \pi_1(BG) = G$$
satisfies
$H_i(BH;\qq) = 0$ for $i \ge 1$.
Then  the higher signatures of $r: M \cup_{\phi} N^- \to X$ and
$s: M \cup_{\psi} N^- \to X$ agree.
\end{corollary}
\proof In view of Lemma \ref{lem: cutting and mapping torus formula in bordims group}
we have to show that the higher signatures of
$(\partial M,\phi^{-1} \circ \psi, r|_{\partial M},K)$ vanish.
The homotopy $K$ yields an element $g \in G$ such that the composition of
$c(g): G \to G \hspace{3mm} g' \mapsto gg'g^{-1}$ with $\pi_1(r|_{\partial M})$ 
agrees with the composition of $\pi_1(r|_{\partial M})$ 
with the automorphism $\pi_1(\phi^{-1} \circ \psi)$. Obviously
$c(g)$ induces an automorphism of $H$. Denote the associated
semi-direct product by $H \rtimes \zz$. There is a group homomorphism
from $H \rtimes \zz$ to $G$ which sends $h \in H$ to $h \in G$ 
and the generator of $\zz$ to $g \in G$. Let $p: H \rtimes \zz \to \zz$ be
the canonical projection. Then  the reference map from the mapping torus 
$T_{\phi^{-1} \circ \psi}$ to $BG$
factorizes as a map $T_{\phi^{-1} \circ \psi} \xrightarrow{t} B(H \rtimes \zz) \to BG$
and the composition $T_{\phi^{-1} \circ \psi} \xrightarrow{t} B(H \rtimes \zz)
\xrightarrow{Bp} B\zz = S^1$ is homotopic to the canonical projection
$\pr: T_{\phi^{-1} \circ \psi} \to S^1$. Notice that $Bp$ induces an isomorphism
$H^*(B\zz;\qq) \to H^*(B(H \rtimes\zz;\qq))$ since $H^i(BH;\qq) = 0$ for $i \ge 1$.
Hence it remains to show that all higher signatures of
$\pr: T_{\phi^{-1} \circ \psi} \to S^1$ vanish. If $1 \in H^0(S^1;\qq) \cong Q$
and $u \in H^1(S^1;\qq) \cong \qq$ are the obvious generators, it remains to prove that
$\sign_1(\pr: T_{\phi^{-1} \circ \psi} \to S^1)$ and
$\sign_u(\pr: T_{\phi^{-1} \circ \psi} \to S^1)$ are trivial. (Recall 
that $\sign_u $ has been defined in the introduction.) These numbers are given by
ordinary signatures of  $T_{\phi^{-1} \circ \psi}$ and $\partial M$.
Since  $T_{\phi^{-1} \circ \psi}$ fibers over $S^1$ and
$\partial M$ is nullbordant, these numbers are trivial. \qed

\begin{example} \label{exa: additivity is true also assumption fails}
\em Consider the situation of Lemma \ref{cor: additivity in the bordism ring}
in the special case, where $\partial M$ looks like 
$\rr\pp^3 \sharp \rr\pp^3 \times N$
for a simply connected oriented closed $(n-3)$-dimensional 
manifold $N$ for $n \ge 6$
such that  $H_{m-1}(N;\qq)$ or $H_{m-3}(N;\cc)$ is non-trivial and 
the map $\pi_1(\partial M) \to \pi_1(M)$ is injective. Recall that
$m$ is the integer satisfying $n = 2m$ resp. $n= 2m+1$.  
The fundamental group of the connected sum
$ \rr\pp^3 \sharp \rr\pp^3$ is the infinite dihedral group 
$D_{\infty} = \zz/2 \ast \zz/2 = \zz \rtimes \zz/2$.
Notice that there is a two-fold covering $S^1 \times S^2 \to \rr\pp^3 \sharp \rr\pp^3$
and the universal covering of $(\rr\pp^3 \sharp \rr\pp^3)^{\sim}$ of
$\rr\pp^3 \sharp \rr\pp^3$ is $\rr \times S^2$. Hence 
$H_{m-1}((\rr\pp^3 \sharp \rr\pp^3 \times N)^{\sim};\cc) 
\cong H_{m-1}(N;\cc) \oplus H_{m-3}(N;\cc)$ is a non-trivial direct
sum of finitely many copies of the trivial $D_{\infty}$-representation $\cc$. Since 
$\alpha_m((\rr\pp^3 \sharp \rr\pp^3 \times N)^{\sim}) =
\alpha_m((S^1 \times S^2 \times N)^{\sim})$ 
\cite[Remark 3.9]{Lott-Lueck (1995)} we conclude
that $\alpha_m((\rr\pp^3 \sharp \rr\pp^3 \times N)^{\sim})$
is different from $\infty^+$
(see \cite[Example 4.3, Theorem 5.4]{Lueck (1997a)} or 
\cite[Theorem 8.7 (9)]{Lueck (2000)}).
Hence Assumption \ref{assumption} is not satisfied because of
Lemma \ref{lem: vanishing middle differential and NS-invariants} and  we cannot 
conclude the additivity of the higher signatures from Corollary 
\ref{cor: Additivity of higher signatures}. But we can conclude
the additivity of the higher signatures from Corollary 
\ref{cor: additivity in the bordism ring} since $H_i(BD_{\infty};\qq) \cong
H_i(B\zz/2);\qq) \oplus H_i(B\zz/2);\qq) = 0$ holds for $i \ge 1$.
\em
\end{example}

More generally one can consider glueing tuples
$(M,r_M,N,r_N,\phi,H)$, which consist of two compact
oriented $n$-dimensional manifolds $M$ and $N$ with boundaries with reference maps
$r_M: M \to X$ and $r_N: N \to X$, an orientation preserving
diffeomorphism $\phi: \partial M \to \partial N$ and a homotopy
$H: \partial M \times [0,1] \rightarrow X$ between $ r_M|_{\partial M} $
 and $r_N|_{\partial N} \circ \phi $, i.e.
$ r_M|_{\partial M} \sim r_N|_{\partial N} \circ \phi$. To such
a glueing tuple one can associate an element
\begin{eqnarray}
\hspace{-8mm} [M,r_M,N,r_N,\phi,H] & := &
[M\cup_{\partial M \times \{0\}}( \partial M \times[0,1])
\cup_{\phi} N^-,r] \hspace{3mm} \in
\Omega_n(X),
\label{bordism element given by glueing tuple}
\end{eqnarray}
where $r$ is constructed from $r_M$, $H$ and $r_N$ in the obvious way.
One gets

\begin{lemma}
\label{lem: additivity for glueing tuples in bordims group}
Let $(M,r_M,N,r_N,\phi,H)$ and $(N,r_N,P,r_P,\psi,G)$ be two
glueing tuples. They can be composed to a glueing tuple
$(M,r_M,P,r_P,\psi\circ\phi,K)$, where $K$ is the composition
of the homotopies $H$ and $G \circ (\phi \times \id)$. Then
we get in $\Omega_n(X)$
\begin{eqnarray*} [M,r_M,P,r_P,\psi\circ\phi,K]
& = & [M,r_M,N,r_N,\phi,H]  + [N,r_N,P,r_P,\psi,G].
\end{eqnarray*}
\end{lemma}
\proof The required bordism has the following underlying manifold.
Take the disjoint union of  $M \times [0,1]$, $N \times [0,1]$ and
$P \times[0,1]$ and identify $(x,t) \in \partial M \times [0,1/3]$ with
$(\phi(x),t) \in N \times [0,1/3]$ and
$(y,t) \in N \times [2/3,1]$ with $(\psi(y),t) \in \partial P \times [2/3,1]$. \qed

\begin{example} \label{exa: non-additivity} \em In odd dimensions additivity
of the higher signatures (sometimes also called the cut and paste property)
fails as bad as possible in the following sense.
Let us consider $m \ge 2$, a finitely presented group $G$
and {\it{any}} element $\omega \in \Omega_{2m+1}(BG)$. Then using the
last (surjective) map of Theorem 3.2  of (\cite{Quinn(1979)}) and
also the isomorphism given at the bottom of page 57 of (\cite{Quinn(1979)})
(or see the Theorem in the appendix by Matthias Kreck),
one can find  a quadruple
$(F,h,r,H)$ for a $2m$-dimensional closed oriented manifold $F$
with reference map $r:F \to BG$
such that $[F,h,r,H] = \omega$ in $\Omega_{2m+1}(BG)$ and $[F,r] = 0$
in $\Omega_{2m}(BG)$ holds. Fix a nullbordism $R: W \to BG$ for $r: F \to BG$.
In the sequel we identity $F = \partial W$. Since $F$ admits a collar
neighborhood in $W$, the inclusion
$F \to W$ is a cofibration and thus  we can extend the homotopy
$H: r \simeq r \circ h$ to a
homotopy $H': R \simeq R'$ for some map $R': W \to BG$ such that
$R'_{|\partial M}= r\circ h$.
Thus we obtain elements $R' \cup_h R : W \cup_h W \to BG$
and $R \cup_{\id} R: W \cup_{\id} W \to BG$ such that $R'\simeq R$. We
conclude from Lemma \ref{lem: cutting and mapping torus formula in bordims group}
\begin{eqnarray}
[R' \cup_{h} R: W \cup_{h} W \to BG] - [R \cup_{\id} R: W \cup_{\id} W \to BG]
& = & \omega.
\label{failure of additivity in odd-dimensions}
\end{eqnarray}

The theorem of Matthias Kreck which he proves in the appendix
shows that for $m \ge 2$ one can arrange
in the situation above that the reference map $r: F \to BG$ is $2$-connected,
provided that $BG$ has finite skeleta.
(Since we only want to have $2$-connected it suffices 
that $BG$ has finite $2$-skeleton.)
Consider the special case $m = 2$ and
$G = \zz$. Choose in \eqref{failure of additivity in odd-dimensions}
the quadruple $(F,h,r,H)$ such that $r: F \to B\zz$ is $2$-connected.
Then $\overline{F}$ is the
universal covering of $F$. We conclude from \cite[Lemma 3.3]{Lott-Lueck (1995)}
that $\alpha_2(\overline{F}) = \alpha_2(\widetilde{S^1}) = \infty^+$.
Therefore Assumption
\ref{assumption} is satisfied for $\overline{\partial W} = \overline{F}$
by Lemma \ref{lem: vanishing middle differential and NS-invariants}.
Notice that there are elements $\omega \in \Omega_5(B\zz)$ whose higher signatures
do not all vanish, for instance $[\cc P^2\times S^1, r]$ where
$r: \cc P^2\times S^1 \rightarrow B \zz = S^1 $ is the projection onto
the second factor. Hence, (for such an example),
 if we set $[M_0,r_0] = [W, R] $ and $[M_1,r_1] = [W', R'] $, then
formula \eqref{failure of additivity in odd-dimensions} and
Theorem  \ref{the: main properties of
sigma(M,r,L)} show that
the right-handside of the formula of Theorem \ref{the: main properties of
sigma(M,r,L)}  \eqref{the: main properties of
sigma(M,r,L): homotopy invariance} is not zero. Thus
 Assumption \ref{assumption} is not enough in the case $n = 2m+1$
(in contrast to the case $n = 2m$ as proven in Theorem
\ref{the: main properties of sigma(M,r)}
\eqref{the: main properties of sigma(M,r): additivity})
to ensure the additivity of the higher signatures.

Counterexamples to additivity in odd dimensions yield also counterexamples
in even dimensions by crossing with $S^1$. In the situation
of \eqref{failure of additivity in odd-dimensions} with
$\omega \in \Omega_{2m+1}(BG)$, we get in
$\Omega_{2m+2}(B(G \times \zz)) = \Omega_{2m+2}(BG \times S^1)$
\begin{eqnarray*}
\lefteqn{[R'\times \id_{S^1} \cup_{h \times \id_{S^1}} R\times \id_{S^1} :
W \times S^1 \cup_{h \times \id_{S^1}} W \times S^1  \to BG \times S^1]} & &
\\
& &  -
[R\times \id_{S^1} \cup_{\id_{\partial W\times S^1}} R\times \id_{S^1} :
W \times S^1 \cup_{\id_{\partial W \times S^1}} W \times S^1  \to BG \times S^1]
\\
& = & \omega \times [\id_{S^1}]
\label{failure of additivity in even-dimensions}
\end{eqnarray*}
and, since the $L-$class of $\omega \times S^1$ may be identified
to the one of $\omega$, for any  $u \in H^*(BG;\qq)$
we have $\sign_{u \times [S^1]}(\omega \times [\id_{S^1}]) = \sign_u(\omega)$
where  $[S^1] \in H^1(S^1;\qq)$ is  the fundamental class. Hence if
$\omega \in \Omega_{2m+1}(BG)$ admits at least a higher signature
which is not zero, then $W\times S^1 \cup_{h\times \id} W\times S^1$
admit a higher signature which is not cut and paste invariant.

\em \end{example}

%%%%%%%%%%%%%%%%%%%%%%%%%%%% Section 2 %%%%%%%%%%%%%%%%%%%%%%%%%%%%%%%%%%%%%%
\typeout{-----------------------  Section 2 ------------------------}

\tit{Computations in symmetric $L$-groups}

In this section we carry out
some algebraic computations and constructions of
classes in symmetric $L$-groups which correspond on the geometric side to
defining higher signatures of manifolds
with boundaries (under Assumption \ref{assumption})
and to glueing processes along boundaries.

We briefly recall some basic facts about (symmetric) Poincar\'e
chain complexes and the (symmetric) $L$-groups
defined in terms of bordism classes of such chain complexes.
For details we refer the reader to \cite{Ranicki(1980a)} and to the
Section $1$ of
\cite{Ranicki(1981)}.

Let $R$ be a ring with involution
$R \to R \hspace{2mm}:\, r \mapsto \overline{r}$.
Two  important examples are the group ring
$\zz G$ with the involution given by
$\overline{g}= g^{-1}$ and the reduced
$C^*$-algebra $C^*_r(G)$ of a group $G$.
Given a left $R$-module $V$, let the
dual $V^*$ be the (left) $R$-module $\hom_R(V,R)$
with the $R$-multiplication given by $(rf)(x) = f(x)\overline{r}$.
Given a chain complex $C_* = (C_*,c_*)$ of (left) R-modules, define
$C^{n-*}$ to be the $R$-chain complex
whose $i$-th chain module is $(C_{n-i})^*$ and
whose $i$-th differential is
$c_{n-i+1}^*: C_{n-i}^* \to C_{n-i+1}^*$.
We call $C_*$  finitely generated projective
if $C_i$ is finitely generated projective for all $i \in \zz$
and vanishes for $i\le 0$.
An $n$-dimensional (finitely generated projective symmetric) Poincar\'e
$R$-chain complex $(C_*,\phi)$ consists
of an $n$-dimensional finitely generated projective $R$-chain complex
$C_*$ together with a $R$-chain homotopy
equivalence $\phi^0_*: C^{n-*} \to C_*$
which the part for $s=0$ of a representative
$\{\phi^s \mid s \ge 0\}$ of an
element in $\phi$ in the hypercohomology group
$Q^n(C_*) = H^n(\zz/2;\hom(C^*,C_*))$. The element $\phi^1$ is
a chain homotopy $(\phi^0)^{n-*} \simeq \phi^0_*$,
where $(\phi^0)^{n-*}$ is obtained
from $\phi^0$ in the obvious way using the
canonical identification $P \to (P^{\ast})^{\ast}$
for a finitely generated projective $R$-module $P$.
The elements $\phi^{s+1}$ are higher
homotopies for $\phi^s_* \simeq (\phi^s)^{n-*}$.

Consider a connected finite $CW$-complex $X$
with universal covering $\widetilde{X}$ and
fundamental group $\pi$. It is an $n$-dimensional Poincar\'e complex
if the (up to $\zz \pi$-chain homotopy well-defined)
$\zz \pi$-chain map $-\cap[X]: C^{n-*}(\widetilde{X})
\to C_*(\widetilde{X})$ is a
$\zz \pi$-chain homotopy equivalence.
Then for any normal covering
$\overline{X} \to X$ with
group of deck transformations $G$, the
fundamental class $[X]$ determines an element
in $\phi \in Q^n(C_*(\overline{X}))$, for which
$\phi^0_*$ is the $\zz G$-chain map induced by $-\cap[X]$ and
$(C_*(\overline{X}),\phi)$ is an $n$-dimensional
Poincar\'e $\zz G$-chain complex
\cite[Proposition 2.1 on page 208]{Ranicki(1980b)}.

The (symmetric) $L$-group $L^n(R)$ is defined
by the algebraic bordism group of $n$-dimensional
finitely generated projective Poincar\'e $R$-chain complexes.
The algebraic bordism relation
mimics the geometric bordism relation. The general philosophy,
which we will frequently use without writing down the details,
is that any geometric construction
for geometric Poincar\'e pairs, such as glueing along a common boundary
with a homotopy equivalence, or  taking mapping tori or writing down
certain bordisms, can be transferred to the category
of algebraic Poincar\'e chain complexes.

However, there is one important
difference between the geometric bordism group $\Omega_n(X)$ and the
$L$-group $L^n(R)$ concerning homotopy invariance. Let $G$ be a group and
let $M, N$ be two  closed oriented
$n$-dimensional manifolds with reference maps
 $r: M \to BG$ and $s: N \to BG$. Suppose that
$f: M \to N$ is a homotopy equivalence such that $s \circ f \simeq r$.
Then this does not imply that the bordism classes $[M,r]$ and
$[N,s]$ agree. But the Poincar\'e $\zz G$-chain complexes $C_*(\overline{M})$ and
$C_*(\overline{N})$ are $\zz G$-chain homotopy equivalent, and
this does imply that their classes in $L^n(\zz G)$ agree
\cite[Proposition 3.2 on page 136]{Ranicki(1980a)}.

The following lemma explains the role of Assumption \ref{assumption}.
Its elementary proof is left to the reader.
\begin{lemma} \label{lem: condition c_m = 0}
Let $C_*$ be a projective $R$-chain complex. Then
the following assertions are equivalent.
\begin{enumerate}

\item $C_*$ is $R$-chain homotopy
equivalent to a $R$-chain complex $D_*$ with trivial $m$-th differential;

\item
$\im(c_m)$ is a direct summand in $C_{m-1}$ where $c_m :\, C_{m-1}
\rightarrow C_{m}$ is the differential;

\item
There is a finitely generated projective $R$-subchain complex $D_* \subset C_*$
with $D_m = \ker(c_m)$, $D_{m-1} \oplus \im(c_{m-1})= C_{m-1}$
and $D_i = C_i$ for $i \not= m,m-1$ such that
the $m$-th differential of $D_*$ is zero and the inclusion $D_* \to C_*$ is
a $R$-chain homotopy equivalence;
\end{enumerate}
\end{lemma}

Fix a non-negative integer $n$. Let $m$ be the integer for which either
$n = 2m$ or $n = 2m+1$.
Next we give an algebraic construction which allows to assign to a
(finitely generated projective symmetric)
Poincar\'e pair $(i_*: C_* \to \overline{C}_*,(\delta\phi, \phi))$
of $n$-dimensional  $R$-chain complexes an element in $L^n(R)$,
provided that $C_*$ is chain homotopy equivalent to a $R$-chain complex
with trivial $m$-th differential.
In geometry this would correspond to assign
to an inclusion $i: \partial M \to M$ of
a manifold $M$ with boundary $\partial M$ together
with a reference map $r: M \to X$
an element in $\Omega_n(X)$, where $C_*$ resp. $\overline{C}_*$ resp. $i_*$ 
plays the role of $C_*(\partial M)$, $C_*(M)$ and $C_*(i)$.
The idea would be to glue some
preferred nullbordism to the boundary. This can be carried out in the
more flexible algebraic setting under rather weak assumptions.

We begin with the case $n = 2m$.
Recall that we assume that $C_*$ is chain homotopy equivalent
to an $R$-chain complex $D_*$
such that $d_m: D_m \to D_{m-1}$ is trivial. Notice that we can arrange that
$D_*$ is $(2m-1)$-dimensional finitely generated projective by
Lemma \ref{lem: condition c_m = 0}.
Fix such a chain homotopy equivalence
$u_*: C_* \to D_*$. Equip $D_*$ with the
Poincar\'e structure $\psi$ induced by
$\phi$ on $C_*$ and $u_*$. Define
$\overline{D_*}$ as the quotient chain complex
of $D_*$ for which $\overline{D}_i = D_i$ if $0 \le i \le m-1$ and
$\overline{D}_i = 0$ otherwise. Let $j_*: D_* \to \overline{D}_*$
be the canonical projection. Notice that it
is indeed a chain map since $d_m$ vanishes.
There is a canonical extension of the Poincar\'e structure
$\psi$ on $D_*$ to a Poincar\'e structure $(\delta\psi,\psi)$
on the pair $j_*: D_* \to \overline{D_*}$, namely, take $\delta\psi$ to be zero.
Now we can glue
the Poincar\'e pairs $(i_*: C_* \to \overline{C}_*,(\delta\phi, \phi))$ and
$(j_*: D_* \to \overline{D_*},(\delta\psi,\psi))$ along the
$R$-chain homotopy equivalence $u_*: C_* \to D_*$
\cite[\S 3]{Ranicki(1980a)}, \cite[1.7]{Ranicki(1981)}.
We obtain a $2m$-dimensional Poincar\'e $R$-chain complex which presents
a class in $L^{2m}(R)$. Since chain
homotopy equivalent Poincar\'e $R$-chain complexes
define the same element in the (symmetric) $L$-groups, this class
is independent of the choice of $u_*: C_* \to D_*$. We denote it by
\begin{eqnarray}
\sigma(i_*: C_* \to \overline{C}_*,(\delta\phi, \phi)) \in L^{2m}(R).
\label{symmetric signature for Poincare pairs for n=2m}
\end{eqnarray}
Notice that a chain homotopy equivalence $u_*: D_* \to E_*$ of $(2m-1)$-dimensional 
chain complexes with trivial $m$-th differential induces a 
chain equivalence $\overline{u}_*: \overline{D}_* \to \overline{E}_*$
such that $u_*$ and $\overline{u}_*$ are compatible with the maps
$D_*\to \overline{D}_*$ and $E_* \to \overline{E}_*$ constructed above. 
Since chain homotopy equivalent Poincar\'e $R$-chain complexes
define the same element in the (symmetric) $L$-groups, the class
defined in \eqref{symmetric signature for Poincare pairs for n=2m}
is independent of the choice of $u_*: C_* \to D_*$.

The proof of the next lemma is straightforward in the sense that one has
to figure out  the argument for the corresponding geometric statements,
what is easy, and then to translate it into the algebraic setting
(see also \cite[Prop. 1.8.2 ii]{Ranicki(1981)}).
\begin{lemma} \label{lem: properties of sigma of a pair for n= 2m}
\begin{enumerate}

\item \label{lem: properties of sigma of a pair for n= 2m: glueing formula}
Let  $(i_*: C_* \to \overline{C}_*,(\delta\phi, \phi))$
and  $(j_*: D_* \to \overline{D}_*,(\delta\psi, \psi))$
be $2m$-dimensional (finitely generated projective symmetric)
Poincar\'e pairs. Let $u_*: C_* \to D_*$ be a $R$-chain equivalence
such that $Q^{2m-1}(u_*): Q^{2m-1}(C_*) \to Q^{2m-1}(D_*)$  maps
$\phi$ to $\psi$. 
Denote by $(E_*,\nu)$ the $2m$-dimensional Poinca\'re chain
complex obtained from $i_*$ and $j_*$ by glueing along $u_*$. 
Suppose that $C_*$ is $R$-chain homotopy equivalent to a
$R$-chain complex with trivial $m$-th differential. Then
we get in $L^{2m}(R)$
\begin{eqnarray*}
\sigma(E_*,\nu) ~ = ~
\sigma(i_*: C_* \to \overline{C}_*,(\delta\phi, \phi))
- \sigma(j_*: D_* \to \overline{D}_*,(\delta\psi, \psi));
\end{eqnarray*}

\item
\label{lem: properties of sigma of a pair for n= 2m: homotopy invariance}
Let  $(i_*: C_* \to \overline{C}_*,(\delta\phi, \phi))$
and  $(j_*: D_* \to \overline{D}_*,(\delta\psi, \psi))$
be two $2m$-dimensional (finitely generated projective symmetric)
Poincar\'e pairs. Let $(\overline{f}_*,f_*): i_* \to j_*$
be a chain homotopy equivalence of pairs, i.e.
$R$-chain homotopy equivalences $\overline{f}_*: \overline{C}_* \to
\overline{D}_*$ and $f_*: C_* \to D_*$ with $\overline{f}_* \circ i_* =
j_* \circ f_*$ such that $Q^n(\overline{f}_*,f_*)$ maps
$(\delta\phi, \phi)$ to $(\delta\psi, \psi)$. Suppose
that $C_*$ is $R$-chain homotopy equivalent to a
$R$-chain complex with trivial $m$-th differential. Then
we get in $L^{2m}(R)$
$$\sigma(i_*: C_* \to \overline{C}_*,(\delta\phi, \phi)) ~ = ~
\sigma (j_*: D_* \to \overline{D}_*,(\delta\psi, \psi)).$$

\end{enumerate}
\end{lemma}

Now the invariant \eqref{sigma(M,r)} is obtained from the invariant 
\eqref{symmetric signature for Poincare pairs for n=2m}
applied to the Poincar\'e pair given by the associated chain complexes.
Theorem \ref{the: main properties of sigma(M,r)}
\eqref{the: main properties of sigma(M,r): glueing formula} and
\eqref{the: main properties of sigma(M,r): homotopy invariance}
follow from
Lemma \ref{lem: properties of sigma of a pair for n= 2m}
\eqref{lem: properties of sigma of a pair for n= 2m: glueing formula} and
\eqref{lem: properties of sigma of a pair for n= 2m: homotopy invariance}.
Theorem \ref{the: main properties of sigma(M,r)}
\eqref{the: main properties of sigma(M,r): additivity}
follows directly from Theorem \ref{the: main properties of sigma(M,r)}
\eqref{the: main properties of sigma(M,r): glueing formula} and
\eqref{the: main properties of sigma(M,r): homotopy invariance} because the
right side of the formula appearing in Theorem \ref{the: main properties of sigma(M,r)}
\eqref{the: main properties of sigma(M,r): glueing formula}
does not involve the glueing diffeomorphism.
Notice that the geometric version of
Lemma \ref{lem: properties of sigma of a pair for n= 2m}
\eqref{lem: properties of sigma of a pair for n= 2m: glueing formula}
has been considered in
Lemma \ref{lem: additivity for glueing tuples in bordims group}

Next we deal with the case $n = 2m+1$.
Recall that we are considering a $(2m+1)$-dimensional finitely
generated  projective Poincar\'e $R$-pair
$(i_*: C_* \to \overline{C}_*,(\delta\phi,\phi))$ and that we assume that
$C_*$ is $R$-chain homotopy equivalent to a $R$-chain complex
$D_*$ with trivial $m$-th differential. Since
$C_* \simeq C^{2m-*} \simeq D^{2m-*}$
holds by Poincar\'e duality, $C_*$ is also $R$-chain homotopy equivalent
to a $R$-chain complex, namely, $D^{2m-*}$
whose $(m+1)$-th differential is trivial.
We conclude from Lemma \ref{lem: condition c_m = 0}
that we can fix a $R$-chain homotopy equivalence $u_*: C_* \to D_*$ to a
$2m$-dimensional finitely generated projective $R$-chain complex $D_*$
such that both $d_{m+1}$ and $d_m$ vanish. This implies
also that $H_m(C_*) \cong H_m(D_*) \cong D_m$ is a
finitely generated projective $R$-module and the Poincar\'e structure
on $C_*$ induces the structure of a $(-1)^m$-symmetric (non-degenerate)
form $\mu$ on $H_m(C_*)$. Recall that a $(-1)^m$-symmetric (non-degenerate)
form $(P,\mu)$ consists of a finitely generated projective $R$-module
$P$ together with an isomorphism $\mu: P \to P^*$ such that the composition
$P \xrightarrow{\cong} (P^*)^* \xrightarrow{\mu^*} P$ of $\mu^*$ with the
canonical isomorphism $P \to (P^*)^*$ is $(-1)^m\cdot \mu$.
The standard $(-1)^m$-symmetric hyperbolic form $H(Q)$ for a
finitely generated projective $R$-module $Q$ is given by
$$\squarematrix{0}{1}{(-1)^m}{0} : H(Q)= Q^* \oplus Q
\to (Q^* \oplus Q)^* = Q \oplus Q^*.$$
A Lagrangian for a $(-1)^m$-symmetric form $(P,\mu)$
is a direct summand $L \subset P$ with inclusion $j: L \to P$
such that the sequence
$0 \to L \xrightarrow{j} P \xrightarrow{j^*\circ \mu} L^* \to 0$ is exact.
Any inclusion $j: L \to P$ of a Lagrangian extends to an isomorphism
of $(-1)^m$-symmetric forms $H(L) \to (P,\mu)$.
A stable Lagrangian for $(P,\mu)$ is a Lagrangian
in $(P,\mu) \oplus H(Q)$ for some
finitely generated projective $R$-module $Q$.
A formation $(P,\mu,K,L)$ consists of a $(-1)^m$-symmetric
(non-degenerate) form $(P,\mu)$ together with two Lagrangians $K,L
\subset P$. A stable formation $(P,\mu,K,L)$ on $(P,\mu)$ is a formation
on $(P,\mu) \oplus H(Q)$ for some finitely generated projective
$R$-module $Q$. For more informations about these notions we refer
to \cite[\S 2]{Ranicki(1980a)}.

There are natural identifications of $L^0(R,(-1)^m)$
with the Witt groups of equivalence classes of $(-1)^m$-symmetric
forms and of $L^1(R,(-1)^m)$ with the Witt group of equivalence classes
of $(-1)^m$-symmetric formations \cite[\S 5]{Ranicki(1980a)}.
There are suspension maps
$L^0(R,(-1)^m) \to L^{2m}(R)$ and
$L^1(R,(-1)^m) \to L^{2m+1}(R)$. These suspension maps
are in contrast to the quadratic $L$-groups not isomorphism
for all rings with involutions, but they are bijective if
$R$ contains $1/2$ \cite[page 152]{Ranicki(1980a)}.
The class of $(C_*,\phi)$ vanishes in $L^{2m}(R)$, an algebraic nullbordism
is given by $(i_* : C_* \to \overline{C}_*,(\delta\phi,\phi))$.
Let $u_*: C_*\to D_*$ be a $R$-chain homotopy equivalence to a $2m$-dimensional
finitely generated projective $R$-chain complex with trivial $m$-th
and $(m-1)$-th differential. Equip $D_*$ with the Poincar\'e structure
$\psi$ induced by the given Poincar\'e structure $\phi$ on
$C_*$ and $u_*$. By doing surgery
on the  projection onto the quotient $R$-chain complex
$D_*|_{m-1}$ whose $i$-th chain module is $D_i$ for $i \le m-1$
and zero otherwise, in the sense of \cite[\S 4]{Ranicki(1980a)}, one sees
that the class of $(C_*,\phi)$ in $L^{2m}(R)$ is the image under
suspension of the element given by the $(-1)^m$-symmetric form
on $H_m(C_*)$.
If $R$ contains $1/2$, the suspension map is bijective.
Hence the $(-1)^m$-symmetric (non-degenerate) form on
$H_m(C_*)$ represents zero
in the Witt group of equivalence classes of $(-1)^m$-symmetric forms.
This shows
\begin{lemma}\label{lem: existence of stable Lagrangian}
Suppose that $(M,r)$ satisfies Assumption \ref{assumption} and that
$1/2 \in R$. Then there exists a stable Lagrangian
$L \subset H_m(C_*(\overline{\partial M}) \otimes_{\zz G} R)$.
\end{lemma}

Now suppose that we have fixed a stable Lagrangian
$L \subset H_m(C_*)$. By adding the $m$-fold suspension of
$H(Q)$ for some finitely generated projective $R$-module $Q$ to
$C_*$, we can arrange that $L \subset H_m(C_*)$ is a (unstable) Lagrangian.
Equip $D_*$ with the Poincar\'e structure $\psi$ induced by $\phi$ and $u_*$.
Let $K \subset H_m(D_*)$ be the Lagrangian given by $L$ and $H_m(u_*)$.
Let $\overline{D}_*$ be the quotient $R$-chain complex
of $D_*$ such that $\overline{D}_i = D_i$ for $i \le m-1$,
$\overline{D}_m = K^*$, $\overline{D}_i = 0$ for $i \ge m+1$,
the $i$-th differential is $d_i: D_i \to D_{i-1}$
for $i \le m-1$ and all other differentials are zero.
Let $j_*: D_* \to \overline{D}_*$ be the  $R$-chain map
which is the identity in dimensions $i \le m-1$ and given by
the obvious composition
$D_m = H_m(D_*) \xrightarrow{\cong} H_m(D^{2m-*}) =  H_m(D_*)^* \to K^*$.
There is a canonical extension of the
Poincar\'e structure $\psi$ on $D_*$
to a structure $(\delta\psi,\psi)$ of a Poincar\'e pair 
on $j_*: D_* \to \overline{D}_*$, namely, put $\delta\psi$ to be zero.
Now we can glue the pairs $i_*: C_* \to \overline{C}_*$ and
$j_*: D_* \to \overline{D}_*$
along $u_*$ to get a $(2m+1)$-dimensional Poincar\'e $R$-chain complex.
Its class in $L^{2m+1}(R)$ does not depend on the choice of $Q$, $D_*$ and $u_*$ and
is denoted by
\begin{eqnarray}
\sigma(i_*: C_* \to \overline{C}_*,(\delta\phi, \phi),L) \in L^{2m+1}(R).
\label{symmetric signature for Poincare pairs for n=2m+1}
\end{eqnarray}

Again the proof of the next lemma is straightforward in the sense that one has
to figure out  the argument for the corresponding geometric statements,
what is easy, and then to translate it into the algebraic setting
(see also \cite[Prop. 1.8.2 ii]{Ranicki(1981)}).
\begin{lemma} \label{lem: properties of sigma of a pair for n= 2m+1}
\begin{enumerate}

\item \label{lem: properties of sigma of a pair for n= 2m+1: glueing formula}
Let  $(i_*: C_* \to \overline{C}_*,(\delta\phi, \phi))$
and  $(j_*: D_* \to \overline{D}_*,(\delta\psi, \psi))$
be $(2m+1)$-dimensional (finitely generated projective symmetric)
Poincar\'e pairs. Let $u_*: C_* \to D_*$ be a $R$-chain equivalences
such that $Q^{2m}(u_*): Q^{2m}(C_*) \to Q^{2m}(D_*)$  maps $\phi$ to $\psi$.
Suppose that $C_*$ is $R$-chain homotopy equivalent to a
$R$-chain complex with trivial $m$-th differential.
Let $K \subset H_m(C_*)$ and $L \subset H_m(D_*)$ be stable Lagrangians
such that $H_m(u_*): H_m(C_*) \to H_m(D_*)$ respects them stably.
Let $(E_*,\nu)$ be the $2m$-dimensional Poinca\'re chain
complex obtained from $i_*$ and $j_*$ by glueing along $u_*$.
Then we get in $L^{2m}(R)$
\begin{eqnarray*}
\sigma(E_*,\nu) ~ = ~
\sigma(i_*: C_* \to \overline{C}_*,(\delta\phi, \phi),K)
- \sigma(j_*: D_* \to \overline{D}_*,(\delta\psi, \psi),L);
\end{eqnarray*}

\item
\label{lem: properties of sigma of a pair for n= 2m+1: homotopy invariance}
Let  $(i_*: C_* \to \overline{C}_*,(\delta\phi, \phi))$
and  $(j_*: D_* \to \overline{D}_*,(\delta\psi, \psi))$
be two $(2m+1)$-dimensional (finitely generated projective symmetric)
Poincar\'e pairs. Let $(\overline{f}_*,f_*): i_* \to j_*$
be a chain homotopy equivalence of pairs, i.e.
$R$-chain homotopy equivalences $\overline{f}_*: \overline{C}_* \to
\overline{D}_*$ and $f_*: C_* \to D_*$ with $\overline{f}_* \circ i_* =
j_* \circ f_*$ such that $Q^{2m+1}(\overline{f}_*,f_*)$ maps
$(\delta\phi, \phi)$ to $(\delta\psi, \psi)$. Suppose
that $C_*$ is $R$-chain homotopy equivalent to a
$R$-chain complex with trivial $m$-th differential.
Let $K \subset H_m(C_*)$ and $L \subset H_m(D_*)$ be stable Lagrangians.
Denote by $K' \subset H_m(D_*)$ the image of $K$ under $H_m(f_*)$.
Then we obtain a stable equivalence class of
formations $(H_m(D_*),\nu,K',L)$. Let
$[H_m(D_*),\nu,K',L] \in L^{2m+1}(R)$ be the image of
the element which is represented
in the Witt group of equivalence classes
of formations under the suspension homomorphism.
Then we get in $L^{2m+1}(R)$
$$\sigma(i_*: C_* \to \overline{C}_*,(\delta\phi, \phi),K) -
\sigma (j_*: D_* \to \overline{D}_*,(\delta\psi, \psi),L)
=[H_m(D_*),\nu,K',L].
 $$

\end{enumerate}
\end{lemma}

Now the invariant \eqref{sigma(M,r,L)} is obtained from the invariant 
\eqref{symmetric signature for Poincare pairs for n=2m+1}
applied to the Poincar\'e pair given by the associated chain complexes.
Theorem \ref{the: main properties of sigma(M,r,L)}
follows from Lemma
\ref{lem: properties of sigma of a pair for n= 2m+1}.

The next example shall illustrate that the choice of the homotopies
$h_M$ and $h_N$ in Theorem \ref{the: main properties of sigma(M,r,L)}
\eqref{the: main properties of sigma(M,r,L): additivity} and of the homotopy $h$ in
Theorem \ref{the: main properties of sigma(M,r,L)}
\eqref{the: main properties of sigma(M,r,L): homotopy invariance} do affect
the terms given by the formations. We are grateful to Michel Hilsum who pointed out to us
that in an earlier version we did not make this point clear enough.

\begin{example} \label{exa: cuting S^1} \em
Put $R= \zz[\zz]$. 
Consider $(D^1,S^0)$ with the following two different reference maps 
$c,e: D^1 = [-1,1] \to B\zz = S^1$, namely $c(s) = \exp(0)$ and $e(s) = \exp(\pi i(s +1))$.
Let $h: D^1 \times [0,1] \to B\zz = S^1$ be the homotopy $e \simeq t$ sending
$(s,t)$ to $\exp(\pi i t(s+1))$. Notice that $c^*E\zz|_{S^0}$ and $t^*E\zz|_{S^0}$
agree and that we can choose therefore for both the same Lagrangian 
$L \subset H_0(\overline{S^0})$. Obviously Assumption
\ref{assumption} is satisfied. We want to show that
$\sigma(\overline{D^1},t,L)$ and $\sigma(\overline{D^1},e,L)$ are not the same elements
in $L^0(\zz[\zz])$. Their difference 
$\sigma(\overline{D^1},t,L) - \sigma(\overline{D^1},e,L)$
is given by  the class of the formation 
$[H_m(C_*(\overline{S^0})),\mu_1, L_0',L_1]$.  From
Theorem \ref{the: main properties of sigma(M,r,L)}
\eqref{the: main properties of sigma(M,r,L): homotopy invariance}
applied to $(f,\partial f) = id: (D^1,S^0) \to (D^1,S^0)$ and $h$,
we can identify the form $(H_0(S^0),\mu)$ with 
$\squarematrix{1}{0}{0}{-1} : \zz[\zz] \oplus \zz[\zz] \to  \zz[\zz] \oplus \zz[\zz]$
and choose $L = \{(x,x) \mid x \in \zz[\zz]\} \subset \zz[\zz] \oplus \zz[\zz]$.
The homotopy $h$ and the identity on $S^0$ induce the $\zz$-automorphism
of $c^*E\zz|_{S^0} = t^*E\zz|_{S^0} = S^0 \times \zz$ which is the identity on
$\{-1\} \times \zz$ and multiplication with $t$ on $\{1\} \times \zz$.
The automorphism of the form $(H_0(S^0),\mu)$ 
induced by $\partial f = \id$ and $h|_{S^0}$ is
$\squarematrix{1}{0}{0}{t}: \zz[\zz] \oplus \zz[\zz] \to  \zz[\zz] \oplus \zz[\zz]$.
Hence the difference $\sigma(\overline{D^1},t,L) - \sigma(\overline{D^1},e,L)$
is represented by the formation $(\zz[\zz] \oplus \zz[\zz],\squarematrix{1}{0}{0}{-1},L',L)$
for $L' = \{(x,tx) \mid x \in \zz[\zz]\} \subset \zz[\zz] \oplus \zz[\zz]$.
But the class of this formation under the isomorphism
$L^0(\zz[\zz])[1/2] \cong L^0(\zz)[1/2] \oplus L^1(\zz)[1/2] \cong \zz[1/2] \oplus 0 = \zz[1/2]$ 
is the generator.

Similarily one can see from 
Theorem \ref{the: main properties of sigma(M,r,L)}
\eqref{the: main properties of sigma(M,r,L): additivity} that
$\sigma(S^1,\id)$ and $\sigma(S^1,c)$ are different for
the reference maps $\id = S^1 \to S^1 = B\zz$ and the constant map
$c: S^1 \to B\zz$ by cutting $S^1$ open along the embedded $S^0 \subset S^1$.
One has to choose homotopies $h_+: i_+ \simeq c : S^1_+ \to B\zz$ and
$h_-: i_- \simeq c : S^1_+ \to B\zz$ for $S^1_{\pm}$ the upper and lower hemispheres
and $i_{\pm}: S^1_{\pm} \to B\zz = S^1$ the inclusion.
Then the term describing $\sigma(\overline{S^1_+ \cup_{\id} S^1_-},\id) -
\sigma(\overline{S^1_+ \cup_{\id} S^1_-},c)$ is given again by a formation which does not 
represent zero in $L^1(\zz[\zz])$. By crossing with $\cc\pp^{2n}$ one gets also examples
in dimensions $4n+1$ of this type because crossing
with $\cc\pp^{2n}$ induces an isomorphism $L^1(\zz[\zz])[1/2] \to 
L^{4n+1}(\zz[\zz])[1/2]$. 
\em
\end{example}

The next lemma is the algebraic version of
Lemma \ref{lem: cutting and mapping torus formula in bordims group}
(see also \cite[Prop. 1.8.2 ii]{Ranicki(1981)}).

\begin{lemma}
\label{lem: properties of sigma of a pair: Additivity and mapping tori}
Let $n$ be any positive integer.
Let  $(i_*: C_* \to \overline{C}_*,(\delta\phi, \phi))$
and  $(j_*: D_* \to \overline{D}_*,(\delta\psi, \psi))$
be two $n$-dimensional (finitely generated projective symmetric)
Poincar\'e pairs. Let $u_*,v_*: C_* \to D_*$ be a $R$-chain equivalences
such that both $Q^n(u_*)$  and $Q^n(v_*)$ map
$(\delta\phi, \phi)$ to $(\delta\psi, \psi)$.
Let $w_*: C_* \to C_*$ be a $R$-chain map with
$u_* \circ w_* \simeq v_*$.
Let $(E_*(u_*),(\delta\nu,\nu)(u_*))$ and
$(E_*(v_*),(\delta\nu,\nu)(v_*))$ respectively
be the $n$-dimensional Poinca\'re chain
complexes obtained from $i_*$ and $j_*$ by glueing along $u_*$ and $v_*$
respectively. Let $(T(w_*), \mu)$ be the algebraic mapping torus
of $w_*$. Its underlying $R$-chain complex is the mapping cone of
$\cone(\id - w_*)$ (cf. \cite[page 264]{Ranicki(1998)}).
Then we get in $L^{n}(R)$
\begin{eqnarray*}
\sigma(E_*(u_*),(\delta\nu,\nu)(u_*)) -
\sigma(E_*(v_*),(\delta\nu,\nu)(v_*)) & = &
\sigma(T(w_*), \mu).
\end{eqnarray*}
\end{lemma}

In general symmetric signatures and
higher signatures are not additive (see Example \ref{exa: non-additivity}).
In the situation of
Lemma \ref{lem: cutting and mapping torus formula in bordims group}
the difference of symmetric signatures (and thus of
higher signatures) is measured by the symmetric signature of
the corresponding mapping torus. If we want to see the difference in
 $L^{n}(C_r^*(G))$,
we only have to consider the algebraic mapping torus as explained
in Lemma
\ref{lem: properties of sigma of a pair: Additivity and mapping tori}.
To detect the image of the class of the mapping torus in $L^{n}(C^*_r(G))$ under
the isomorphism $\sign: L^{n}(C^*_r(G)) \to K_0(C^*_r(G))$ the formula
\cite[Proposition 4.3]{Ranicki(1980a)} is useful. It reduces the computation of
the difference of the element
$[r: M \cup_{\phi} N^- \to X] - [s: M \cup_{\psi} N^- \to X]$ under the composition
$\Omega_{n}(BG) \xrightarrow{D} K_n(BG) \xrightarrow{A} K_n(C^*_r(G))$
to an expression which only involves
the chain complex of $C_*(\overline{\partial M})$ and the map induced by
the automorphism  $\phi^{-1} \circ \psi$ in a rather close range
around the middle dimension.

\begin{remark} \label{rem: vanishing of all higher signatures} \em
Let $Z$ be a closed oriented $n$-dimensional manifold with a reference map
$r: Z \to BG$.
Suppose that we have for the $m$-th Novikov-Shubin invariant
$\alpha_m(\overline{Z}) = \infty^+$
in the case $n = 2m-1$ and, in the case $n = 2m$, 
we have $\alpha_m(\overline{Z}) = \infty^+$ and
for the $m$-th $L^2$-Betti number
$b_m^{(2)}(\overline{Z}) = 0$. Then we conclude from the arguments above
and Lemma \ref{lem: vanishing middle differential and NS-invariants}
that $\sigma: \Omega_n(BG) \to L^n(C_r^*(G;\rr))$ maps $[Z,r]$ to zero.
Namely, we have constructed an explicit algebraic nullbordism above. 
Hence we conclude that all higher signatures of $[Z,r]$ vanish if the assembly map
$A_{\rr}: KO_n(BG) \otimes_{\zz} \qq \to KO_n(C^*_r(G;\rr))$ is injective. 
This follows from the discussion in the introduction.
\em
\end{remark}

%%%%%%%%%%%%%%%%%%%%%%%%%%%% Section 3 %%%%%%%%%%%%%%%%%%%%%%%%%%%%%%%%%%%%%%
\typeout{-----------------------  Section 3 ------------------------}

\tit{Novikov-Shubin invariants}

Next we reformulate (following \cite{Leichtnam-Lott-Piazza(1999)}) 
the condition that the middle differential vanishes
in terms of spectral invariants.

Let $\caln(G)$ be the von Neumann algebra associated to $G$.
Let $M$ be a closed Riemannian manifold with normal covering $\overline{M} \to M$
with deck transformation group $G$. Let $\nu$ be the flat
$C^*_r(G)$-bundle over $M$ whose total space is
$\overline{M} \times_G C^*_r(G)$. Let
$H^m(M;\nu)$ and $\overline{H}^m(M;\nu)$ resp.
be the unreduced and reduced
$m$-th cohomology of $M$, i.e. $\ker(d^m)/\im(d^{m-1})$ and
$\ker(d^m)/\overline{\im(d^{m-1})}$ resp. for $d$ the differential in the
deRham complex $\Omega^*(M;\nu)$ of Hilbert $C^*_r(G)$-modules.
The next lemma is contained in Lemmas 2.1 and 2.3 of
\cite{Leichtnam-Lott-Piazza(1999)}.

\begin{lemma}
\label{lem: vanishing middle differential and NS-invariants}
The following assertions are equivalent
for an integer $m$.

\begin{enumerate}

\item \label{lem: vanishing middle differential and NS-invariants: C^*-bundle}
The canonical projection
$H^m(\overline{M};\nu) \to \overline{H}^m(\overline{M};\nu)$
is bijective;

\item \label{lem: vanishing middle differential and NS-invariants: C^*}
The $C^*_r(G)$-chain complex $C_*(\overline{M}) \otimes_{\zz G} C^*_r(G)$
is $C^*_r(G)$-chain homotopy equivalent to a finitely generated projective
$C^*_r(G)$-chain complex
$D_*$ whose $m$-th differential $d_m: D_m \to D_{m-1}$ is trivial;

\item \label{lem: vanishing middle differential and NS-invariants: caln(G)}
The $\caln(G)$-chain complex $C_*(\overline{M}) \otimes_{\zz G} \caln(G)$
is $\caln(G)$-chain homotopy equivalent to a finitely generated projective
$\caln(G)$-chain complex
$D_*$ whose $m$-th differential $d_m: D_m \to D_{m-1}$ is trivial;

\item \label{lem: vanishing middle differential and NS-invariants: NS}
The Novikov Shubin invariant $\alpha_m(\overline{M})$
is $\infty^+$ (see \cite{Lott-Lueck (1995)});

\item \label{lem: vanishing middle differential and NS-invariants: Laplacian}
The Laplacian acting on $L^2(\overline{M}, \Omega^{m-1})/\ker(d^{m-1})$
has a strictly positive spectrum.
\end{enumerate}
\end{lemma}
\proof
\eqref{lem: vanishing middle differential and NS-invariants: C^*-bundle}
$ \Leftrightarrow $
\eqref{lem: vanishing middle differential and NS-invariants: C^*}
We can interprete the (a priori purely algebraic) $C^*_r(G)$-cochain complex
$\hom_{\zz G}(C_*(\overline{M}),C^*_r(G)$ as cochain complexes
of Hilbert $C^*_r(G)$-chain complexes with adjointable morphisms as differentials
by the identification of each cochain module
with the direct sum of finitely
many copies of $C_*^r(G)$ using cellular $\zz G$-basis.
There is a $C^*_r(G)$-chain homotopy equivalence
(by bounded chain maps and homotopies)
$\Omega^*(M;\nu) \to \hom_{\zz G}(C_*(\overline{M}),C^*_r(G))$.
Hence the image of the $(m-1)$-th differential in
$\Omega^*(M;\nu)$ is closed if and only if the same is true for
the one in $\hom_{\zz G}(C_*(\overline{M}),C^*_r(G))$.
The image of a differential in $\hom_{\zz G}(C_*(\overline{M}),C^*_r(G))$
is closed if and only if the image is a direct
summand in the purely algebraic
sense \cite[Corollary 15.3.9]{Wegge-Olsen(1993)}.
But this is equivalent to the assertion that
$\hom_{\zz G}(C_*(\overline{M}),C^*_r(G))$ is $C^*_r(G)$-chain homotopy
equivalent to a finitely generated projective $C^*_r(G)$-cochain complex  whose
$(m-1)$-th codifferential is trivial by
Lemma \ref{lem: condition c_m = 0}. This is true if and only
if $C_*(\overline{M}) \otimes_{\zz G} C^*_r(G)$ is $C^*_r(G)$-chain homotopy
equivalent to $C^*_r(G)$-chain complex with trivial $m$-th differential.
\\[2mm]
\eqref{lem: vanishing middle differential and NS-invariants: C^*}
$ \Rightarrow $
\eqref{lem: vanishing middle differential and NS-invariants: caln(G)}
is obvious.
\\[2mm]
\eqref{lem: vanishing middle differential and NS-invariants: caln(G)}
$\Leftrightarrow $
\eqref{lem: vanishing middle differential and NS-invariants: NS}
follows directly from the interpretation of Novikov-Shubin invariants in
terms of the homology of
$C_*(\overline{M}) \otimes_{\zz G} \caln(G)$
\cite{Lueck (1997a)}.
\\[2mm]
\eqref{lem: vanishing middle differential and NS-invariants: NS}
$\Leftrightarrow $
\eqref{lem: vanishing middle differential and NS-invariants: Laplacian}
follows from the fact that the dilatational equivalence
class of the spectral density function of the simplicial
$m$-th codifferential and the analytic $m$-th codifferential
agree. \cite{Gromov-Shubin(1991)}.\\[2mm]
\eqref{lem: vanishing middle differential and NS-invariants: Laplacian}
$\Leftrightarrow $
\eqref{lem: vanishing middle differential and NS-invariants: C^*-bundle}
Assertion \eqref{lem: vanishing middle differential and NS-invariants: C^*-bundle}
can be reformulated to the statement that the spectrum of $(d^{m-1})^*d^{m-1}$
for $d^{m-1}$ $(m-1)$-th differential in
$\Omega^*(M;\nu)$ has a gap  at zero. 
But this spectrum is the same as the
spectrum of $(d^{m-1})^*d^{m-1}$ for $d^{m-1}$ the $(m-1)$-th differential
in the deRham complex
$L^2\Omega^*(\overline{M})$ of Hilbert spaces which has a gap at zero
if and only if \eqref{lem: vanishing middle differential and NS-invariants: Laplacian} 
is true. \qed

%%%%%%%%%%%%%%%%%%%%%%%%%%%% References  %%%%%%%%%%%%%%%%%%%%%%%%%%%%%%%
\typeout{-----------------------  References ------------------------}

%%%%%%%%%%%%%%%%%%%%%%%%%%%%%% Appendix %%%%%%%%%%%%%%%%%%%%%%%%%

\typeout{------------------------  Appendix ---------------}
\vspace{5mm}
\begin{center}\begin{Large}\begin{bf}
Appendix: Mapping tori of special diffeomorphisms\\
by\\
Matthias Kreck
\end{bf}\end{Large}\end{center}

In this appendix we consider the image of the bordism group of
diffeomorphisms on smooth manifolds over a $CW$-complex $X$
with finite skeleta under the mapping torus construction.
By a diffeomorphism over $X$ we mean a quadrupel $(M, f, g, h)$,
where $M$ is a closed oriented smooth manifold,
$g$ an orientation preserving diffeomorphism on $M$,
$f: M \to X$ a continuous map and $h$ a homotopy between $f$ and $f\circ g$.
The role of the homotopy $h$ becomes clear if we consider
the mapping torus $M_g := M \times [0,1]/_{(g(x),0) \sim (x,1)}$,
which is by projection to the second factor a smooth fibre
bundle over $S^1$. Then $h$ allows an extension of $f$ on the
fibre over $0$ to a map $\bar h([(x,t)] ) := h(x,t)$ and any
such extension gives a homotopy $h$ with the properties above.

Following \cite{K1} we denote the bordism group of these
quadrupels by $\Delta _n(X)$. Let $\Omega _{n+1}(X)$ 
be the bordism group of  oriented smooth manifolds with reference 
map to $X$. The mapping torus construction
above gives a homomorphism $ \Delta _{n}(X) \to \Omega _{n+1}(X)$.
It was shown in \cite{K1} for $X$ simply connected and in \cite{Q1} for
general $X$ that for $n$ even this map is surjective.
Recently Wolfgang L\"uck and Eric Leichtnam asked whether
the same statement holds if we only allow quadrupels
where the map $f$ is $2$-connected
and $(M,f)$ represents zero in $\Omega_n(X)$.
We call such a quadrupel a \emph{special diffeomorphism} over
$X$ and the subset of $\Delta _n(X)$ represented by special
diffeorphisms by $S \Delta _n(X)$
(it is not clear to the author whether this subset forms a subgroup).
For $X$ simply connected one can conclude from \cite[\S 9 ]{K1},
that $S \Delta _{2n}(X) \to \Omega _{2n+1}(X)$ is surjective.
In this note we generalize this to arbitrary complexes $X$.
\\[5mm]
{\bf Theorem:} {\em Let $X$ be a $CW$-complex with finite skeleta.
For $n \ge 2$ the mapping torus construction gives a
surjection $S \Delta _{2n}(X) \to \Omega _{2n+1}(X)$.
}\\[4mm]
{\bf Proof:} Let $(N,g)$ be an element of $\Omega _{2n+1}(X)$. Consider a
representative $(M,r)$ of $0$ in $\Omega _{2n}(X)$. We use the
language and results from \cite{K2}. We consider the
fibration $p_2 : X \times BSO \to BSO$ and denote it by $B$.
The map $r \times \nu : M \to B$, where $\nu$ is the normal Gauss map,
is a normal $B$-structure. By \cite[Corollary 1]{K2} we can
replace $(M,r\times \nu)$ up to bordism by a
$n$-equivalence $r' \times \nu ': M'\to X \times BSO$
giving a normal $(n-1)$-structure on $M'$. In particular
$r' : M' \to X$ is $2$-connected.

Now we form the disjoint union $(M' \times I) + N$ and
consider the map $q: (M' \times I) + N \to X$ given by $r' p_1$ and $g$,
where $p_1: M' \times I \to M'$ is the projection. We
want to replace this manifold by a manifold $W$
diffeomorphic to $M \times I$ which is bordant  relative boundary over $X$ to
$(M'\times I) + N$. If this is possible we are finished
since then we glue the two boundary components of $W$ and the
maps together to obtain a mapping torus and a map to $X$.
This is bordant over $X$ to $(N,g)$ since it is bordant
to $((M' \times S^1) + N,r' p_1 + g)$
(note that $(M' \times S^1, r'p_1)$ is zero bordant over $X$).

This idea does not work directly. What we will prove is that
there is a bordism $W$ between $M' \# m(S ^n \times S^n)$ and
$M' \# m(S ^n \times S^n)$ for some $m$ equipped with a map to $X$
which on the two boundary components is the composition of
the projection from $M \# m(S^n \times S^n)$ to $M$ and $r'$,
such that $W$ is diffeomorphic to $(M' \# m(S^n \times S^n)) \times I$.
We further achieve that the manifold obtained by glueing the
boundary components of $W$ together is over $X$ bordant to $(N,g)$.
This is by the considerations above enough to prove the theorem,
since our map from $M \# m(S ^n \times S^n)$ to $X \times BSO$ is
again a $n$-equivalence.

That this indirect way works follows from \cite[Theorem 2]{K2},
which says that we can replace $(M' \times I) + N$
by a sequence of surgeries over $X \times BSO$ and
compatible subtractions of tori by an s-cobordism
$W$ between $M' \# m(S^n \times S^n)$ and
$M' \# m(S^n \times S^n)$ (the fact that the number
of $S^n \times S^n$'s one has to add by Theorem 2
to the boundary components of $W$ is equal follows
from the equality of the Euler characteristic of the
two boundary components). If $n >2$ the $s$-cobordism
theorem implies $W$ diffeomorphic to
$(M' \# m(S^n \times S^n)) \times I$. If $n =2$ the same
is true by the stable $s$-cobordism theorem of \cite{Q2}
after further stabilization of $W$ by forming $k$ times
a "connected sum" between $(S^2 \times S^2) \times I$
and $W$ along an embeded arc joining the two
boundary components of $W$. To finish the argument
one has to note from the definition of compatible subtraction
of tori that this process does not affect the bordism class
over $X$ for the manifold obtained by glueing the
two boundary components together.

To see this we recall the definition of subtraction of tori.
Consider two disjoint embeddings of $S^n \times D^{n+1}$
into $W$ such that the map to $X$ is constant on
both $S^n \times 0$'s. Join each of these embedded tori
by an embedded $I \times D^{2n}$ with the two boundary
components and subtract the interior of these embedded
submanifolds to obtain $W'$. This is the subtraction of a pair
of tori used in \cite[Theorem 2]{K2}. The boundary of $W$ consists
of two copies of $M \# (S^n \times S^n)$. There is an obvious
bordism over $X$ between the manifold obtained from $W$ by
identifying the two boundary components and the manifold
obtained from $W'$ by identifying the two boundary components.\qed

{\bf Remark:}
In general it is difficult to say much about the special diffeomorphism
whose mapping torus is bordant to the given pair $(N,g)$.
The main difficulty is the determination of the diffeomorphism.
One can obtain some information on $M'$. For example
if $X = S^1$ and $n=2$ the proof above shows that we can take
for $M'$ the following manifold:
$S^1 \times S^3 \# \mathbb {CP}^2 \# \bar {\mathbb {CP}}^2$ and
thus the special diffeomorphism lives on
$S^1 \times S^3 \# \mathbb {CP}^2 \# \bar {\mathbb {CP}}^2 \# m(S^2 \times S^2)$
for some unknown integer $m$. More generally in dimension $4$
for an arbitrary $X$ one can use instead of $S^1 \times S^3$
the boundary of any thickening of the $2$-skeleton of $X$ in
$\mathbb R ^5$.

\end{document}